\def\marker{\>\hbox{${\vcenter{\vbox{
    \hrule height 0.4pt\hbox{\vrule width 0.4pt height 6pt
    \kern6pt\vrule width 0.4pt}\hrule height 0.4pt}}}$}\>}
\def\gpic#1{#1
     \medskip\par\noindent{\centerline{\box\graph}} \medskip}
\documentclass[12pt]{article}
\usepackage{amsthm}
\usepackage{url}
\usepackage{verbatim}
\newtheorem*{question}{Question}
\newtheorem{theorem}{Theorem}
\newtheorem{lemma}[theorem]{Lemma}

\newtheorem{corollary}[theorem]{Corollary}
\newtheorem{conjecture}[theorem]{Conjecture}
\newtheorem{definition}[theorem]{Definition}
\newtheorem{example}[theorem]{Example}

\def\st{\colon\,}
\def\VEC#1#2#3{#1_{#2},\ldots,#1_{#3}}
\def\CL#1{\left\lceil #1\right\rceil}
\def\esub{\subseteq}

\begin{document}

\title{Classes of $3$-regular graphs that are $(7,2)$-edge-choosable}
\author{Daniel W. Cranston\thanks{DIMACS, Rutgers University, Piscataway, NJ
         (\texttt{dcransto@dimacs.rutgers.edu}).}
   \and Douglas B. West\thanks{University of Illinois, Urbana, IL
         (\texttt{west@math.uiuc.edu}).   Research partially supported by
         the National Security Agency under Award No.~H98230-06-1-0065.}} 
\maketitle

\begin{abstract}
A graph is $(7,2)$-edge-choosable if, for every assignment of lists of size $7$
to the edges, it is possible to choose two colors for each edge from its list
so that no color is chosen for two incident edges.  We show that every
$3$-edge-colorable graph is $(7,2)$-edge-choosable and also that many
non-$3$-edge-colorable $3$-regular graphs are $(7,2)$-edge-choosable.
\end{abstract}

\section{Introduction}

An \textit{$s$-tuple coloring} of a graph is an assignment of $s$-sets (of
colors) to the vertices such that adjacent vertices receive disjoint sets.
This notion was introduced by Stahl~\cite{stahl}.  A graph $G$ is
\textit{$(r,s)$-colorable} if it has an $s$-tuple coloring using altogether at
most $r$ colors, and the \textit{$s$-set chromatic number} is the least $r$
such that $G$ is $(r,s)$-colorable.  The \textit{fractional chromatic number}
$\chi^*(G)$ (also called \textit{set chromatic number},
\textit{multicoloring number}, or \textit{ultimate chromatic number}) is
$\inf\{r/s\st G \textrm{ is $(r,s)$-colorable}\}$.  Since the ordinary
chromatic number is achievable as such a ratio (with $s=1$), always
$\chi^*(G)\le \chi(G)$.

The parameter $\chi^*$ is worthy of study for various reasons.  It is the
linear programming relaxation of the ordinary chromatic number of a graph (and
hence the infimum is a minimum); many papers have studied it from this
viewpoint.  Furthermore, Klostermeyer and Zhang~\cite{klost} observed that
$\chi^*(G)$ is the minimum ratio $r/s$ such that there is a homomorphism from
$G$ to the Kneser graph $K(r,s)$, where $K(r,s)$ is the graph on the
$s$-subsets of an $r$-element set in which adjacency is defined by disjointness.

A modern variation for coloring problems restricts the colors available at
vertices.  A \textit{list assignment} $L$ gives each vertex $v$ a list $L(v)$
of available colors.  A good coloring must satisfy the additional requirement
that any colors used at $v$ must lie in $L(v)$.  A graph $G$ is
\textit{$k$-choosable} if an ordinary proper coloring can be chosen from the
lists whenever the list have size at least $k$, and the \textit{choosability}
(or \textit{list chromatic number}) $\chi_\ell(G)$ is the least $k$ such that
$G$ is $k$-choosable.  More generally, $G$ is \textit{$s$-set $k$-choosable} if
an $s$-set coloring can be chosen from any assignment of lists of size $k$, and
the \textit{$s$-set choosability} is the minimum $k$ such that $G$ is
$s$-set $k$-choosable.

When coloring edges, the sets chosen for incident edges must be disjoint,
so the edge-coloring problems reduce to the vertex coloring problems mentioned
above on the line graph.  In particular, in this paper we study the
\textit{$2$-set edge-choosability} of $3$-regular graphs.  Since extra colors
can be discarded without making the choosing of a proper coloring easier, it
suffices to study \textit{$r$-uniform} list assignments, where each list has
size $r$.  Given a list assignment $L$ on the edges of $G$, a selection of $s$
colors from each list so that incident edges have disjoint lists is an
\textit{$s$-set $L$-coloring} of $E(G)$.

\begin{question}
What is the least $r$ such that, for every $3$-regular graph $G$, every
$r$-uniform list assignment $L$ on the edges of $G$ admits a $2$-set
$L$-coloring?
\end{question}

That is, we seek the least $r$ such that every $3$-regular graph is
$(r,2)$-edge-choosable.  On his website, in a ``Problem of the Month'', Bojan
Mohar~\cite{mohar} asked for this value.  He conjectured that every $3$-regular
graph is $(7,2)$-edge-choosable.  

A generalization of Brooks' Theorem implies that every $3$-regular graph is
$(8,2)$-edge-choosable.  Tuza and Voigt~\cite{tuza2} proved that if a connected
graph $G$ is not complete and is not an odd cycle, then $G$ is
$(\Delta(G)m,m)$-choosable whenever $m\geq 1$.  Since the line graph of a
$3$-regular graph has maximum degree $4$, every $3$-regular graph is thus
$(8,2)$-edge-choosable.  

On the other hand, it is also easy to construct a $3$-regular graph that is
not $(6,2)$-edge-choosable (in fact, not even $(6,2)$-edge-colorable).
It is well known that the smallest $3$-regular graph $G$ that is not
$3$-edge-colorable is formed from two copies of $K_4$ by subdividing one
edge in each and making the two new vertices adjacent.  A $(6,2)$-edge-coloring
of $G$ would put a total of $30$ colors on the $15$ edges.  With only six
colors available, each color class would have to be a perfect matching.
Since every perfect matching in $G$ uses the central cut-edge, this contradicts
that we choose only two colors on it. 

Thus Mohar's conjecture on the $2$-set choosability of $3$-regular graphs
is sharp if true.

\bigskip
In this paper, we show that every $3$-edge-colorable graph is
$(7,2)$-edge-choosable and that many $3$-regular graphs that are not
$3$-edge-colorable are also $(7,2)$-edge-choosable.  

Some classes of $3$-regular graphs are known to be $(6,2)$-edge-choosable.
Ellingham and Goddyn~\cite{ellingham} showed that planar $d$-regular
$d$-edge-colorable multigraphs are $d$-edge-choosable; by doubling the edges,
this implies that planar $3$-regular graphs are $(6,2)$-edge-choosable.
Haxell and Naserasr~\cite{haxell} showed that the Petersen graph is
$(6,2)$-edge-choosable.   Both of these results use the Alon-Tarsi Theorem and
thus provide only existence proofs.  Our proofs of $(7,2)$-edge-choosability
provide a fairly simple algorithm for choosing a $2$-set edge-coloring from
lists of size $7$ on the edges.

An edge $e$ of a $3$-regular graph is incident to four other edges.  Choosing
colors for them could forbid eight colors from usage on $e$.  Our main idea is
to show that we can choose the four colors on two of these incident edges using
only three colors from $L(e)$.  Any time we reduce the number of available
colors at $e$ by less than the number of colors we choose on edges incident to
$e$, we say that we have \textit{saved a color on $e$}.  In particular, when we
save a color on $e$ while choosing four colors on two edges incident to $e$, we
retain a list of (at least) four colors available for $e$.  We will apply this
repeatedly to choose colors on some edges so that the remaining graph is
$(4,2)$-edge-choosable and retains lists of size $4$.

\section{The Key Idea}
Our first lemma is a simple form of our main tool.  It is a generalization of
the well-known result~\cite{tuza} that even cycles are $(2m,m)$-edge-choosable.
To understand the proof, it may be helpful to picture the case when $A\cup B$
forms an even cycle.  In general, however, $B$ need not be a matching.
We will be saving a color on edges of $B$.
\begin{lemma}
\label{key lemma}
Let $A$ and $B$ be sets of $k$ edges in a graph, with
$A = \{a_1, \ldots, a_k\}$ and $B = \{b_1, \ldots, b_k\}$.  Suppose that $A$
is a matching and that $b_i$ is incident to $a_i$ and $a_{i+1}$ but not
to any other edge in $A$ (the indices are viewed modulo $k$).  From a uniform
list assignment $L$ on the edges, one can choose one color at each edge of $A$
so that for each $i$, together $a_i$ and $a_{i+1}$ receive at most one color
from $L(b_i)$.
\end{lemma}

\begin{proof}
We choose a color $\phi(e)$ for each edge $e$ in $A$.  If the lists for all
edges in $A\cup B$ are identical, then use the same color on each edge of $A$.  
If they are not all identical, then they differ for two consecutive edges in the
cyclic list $a_1, b_1, a_2, b_2, \ldots a_k, b_k$.  We may index the edges so
that these are $a_1$ and $b_k$.  

Since the lists have the same size, we may choose $\phi(a_1)\in L(a_1)-L(b_k)$.
If $\phi(a_1)\notin L(b_1)$, then choose $\phi(a_2)$ from $L(a_2)$ arbitrarily.
If $\phi(a_1)\in L(b_1)\cap L(a_2)$, then let $\phi(a_2) = \phi(a_1)$.  
Finally, if $\phi(a_1)\in L(b_1) - L(a_2)$, then choose
$\phi(a_2)\in L(a_2)-L(b_1)$.  In each case, at most one of the colors chosen
for the edges $a_1$ and $a_2$ incident to $b_1$ is in $L(b_1)$.

Continue in the same manner choosing colors for edges $a_3, \ldots a_k$ so that
at most one color from $L(b_i)$ is used on $a_i$ and $a_{i+1}$.  At the end,
at most one color from $L(b_k)$ appears in $\{\phi(a_k),\phi(a_1)\}$, since
$\phi(a_1)\notin L(b_k)$.
\end{proof}

\begin{corollary}
\label{even cycles}
{\rm(Tuza and Voigt \cite{tuza})}
Even cycles are $(2m,m)$-edge-choosable.
\end{corollary}
\begin{proof}
Partition the edge set into two matchings, $A$ and $B$.  Choose one color for
each edge of $A$ as guaranteed by Lemma~\ref{key lemma}.  Repeat this step $m$
times.  (With each repetition, we may need to discard colors from some lists so
that the remaining lists have equal sizes.)  Now each edge of $B$ has at least
$m$ remaining available colors, which we use on those edges.
\end{proof}


It is not immediately obvious that Lemma~\ref{key lemma} implies anything more
than Corollary~\ref{even cycles}.  Its power lies in carefully choosing the
edge sets $A$ and $B$, as we show in the next proof.  Theorem~\ref{3ec} is a
special case of Theorem~\ref{mainthm}, but it is useful to prove
Theorem~\ref{3ec} independently because it is much simpler than the general
result and yet illustrates our main technique.

\begin{theorem}\label{class1}
\label{3ec}
Every $3$-edge-colorable graph is $(7,2)$-edge-choosable.
\end{theorem}
\begin{proof}
Let $G$ be a 3-edge-colorable graph.  We can raise the degree of each vertex
that has degree less than $3$ by taking two disjoint copies of $G$ and making
the two copies of each deficient vertex adjacent.  A proper $3$-edge-coloring of
$G$ extends to a proper $3$-edge-coloring of the new graph.  If $G$ has 
minimum degree $k$, then we obtain a $3$-regular $3$-edge-colorable supergraph
of $G$ after $3-k$ iterations of this transformation.  We may thus assume that
$G$ is $3$-regular, since every subgraph of a $(7,2)$-edge-choosable graph is
$(7,2)$-edge-choosable.

Let $J$, $K$, and $L$ be the three color classes in a proper $3$-edge-coloring
of $G$.  Since $G$ is $3$-regular, $|J|=|K|=|L|$.  We apply
Lemma~\ref{key lemma} in two phases.  The first phase uses even cycles in
$J\cup K$, with the portion in $J$ as $A$ and the portion in $K$ as $B$.
It chooses one color for each edge of $J$.  The lists of colors remaining
available have size at least 6 for edges of $J$ and $K$ and size at least 5 for
edges of $L$.

The second phase applies Lemma~\ref{key lemma} using even cycles in $J\cup L$,
with the portion in $J$ as $A$ and the portion in $L$ as $B$.  Since the
lemma requires equal list sizes, discard colors from the lists on edges of $J$
and $L$ to reduce to size 5 before applying the lemma.  We have now chosen two
colors for each edge in $J$ but none for any edge in $K$ or $L$.  
There remain at least four colors available at each edge of $K\cup L$.

Since the edges of $K\cup L$ form vertex-disjoint even cycles, applying
Lemma~\ref{even cycles} to each cycle completes the desired coloring.
\end{proof}

Like the iteration in Corollary~\ref{even cycles}, iterating the argument of
Theorem~\ref{class1} shows that every $3$-edge-colorable graph is
$(\CL{7m/2},m)$-edge-choosable, for each positive integer $m$.

\section{The Main Result}
The ideas in Theorem~\ref{3ec} can be applied to prove that many other graphs
are (7,2)-edge-choosable.  Some of these are ``snarks''.  A \textit{snark} is a
$2$-edge-connected $3$-regular graph with edge-chromatic number $4$, girth at
least $5$, and cyclic-connectivity at least $4$.
The smallest snark is the Petersen graph.  The drawing of it in Figure 1 has
exactly one vertical edge $c$ and two horizontal edges labeled $c'$ and $b$.
Initial interest in snarks was due to the equivalence of the Four Color Theorem
with the statement that no planar snarks exist.  They remain interesting
because many important conjectures reduce to the special case of snarks.

\begin{figure}[htb]
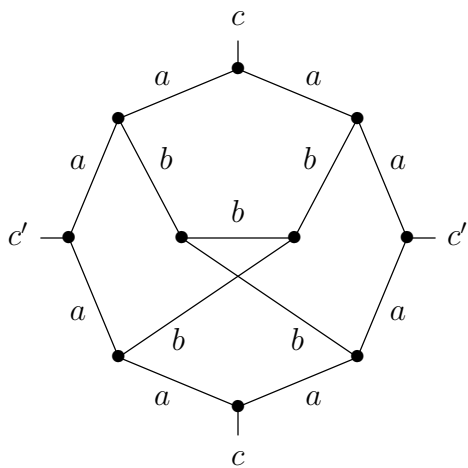

\vspace{-.3 in}
\gpic{
\expandafter\ifx\csname graph\endcsname\relax \csname newbox\endcsname\graph\fi
\expandafter\ifx\csname graphtemp\endcsname\relax \csname newdimen\endcsname\graphtemp\fi
\setbox\graph=\vtop{\vskip 0pt\hbox{%
    \graphtemp=.5ex\advance\graphtemp by 0.265in
    \rlap{\kern 1.150in\lower\graphtemp\hbox to 0pt{\hss $\bullet$\hss}}%
    \graphtemp=.5ex\advance\graphtemp by 0.525in
    \rlap{\kern 1.775in\lower\graphtemp\hbox to 0pt{\hss $\bullet$\hss}}%
    \graphtemp=.5ex\advance\graphtemp by 1.150in
    \rlap{\kern 2.035in\lower\graphtemp\hbox to 0pt{\hss $\bullet$\hss}}%
    \graphtemp=.5ex\advance\graphtemp by 1.775in
    \rlap{\kern 1.775in\lower\graphtemp\hbox to 0pt{\hss $\bullet$\hss}}%
    \graphtemp=.5ex\advance\graphtemp by 2.035in
    \rlap{\kern 1.150in\lower\graphtemp\hbox to 0pt{\hss $\bullet$\hss}}%
    \graphtemp=.5ex\advance\graphtemp by 1.775in
    \rlap{\kern 0.525in\lower\graphtemp\hbox to 0pt{\hss $\bullet$\hss}}%
    \graphtemp=.5ex\advance\graphtemp by 1.150in
    \rlap{\kern 0.265in\lower\graphtemp\hbox to 0pt{\hss $\bullet$\hss}}%
    \graphtemp=.5ex\advance\graphtemp by 0.525in
    \rlap{\kern 0.525in\lower\graphtemp\hbox to 0pt{\hss $\bullet$\hss}}%
    \special{pn 8}%
    \special{pa 1150 265}%
    \special{pa 1775 525}%
    \special{pa 2035 1150}%
    \special{pa 1775 1775}%
    \special{pa 1150 2035}%
    \special{pa 525 1775}%
    \special{pa 265 1150}%
    \special{pa 525 525}%
    \special{pa 1150 265}%
    \special{fp}%
    \special{pa 1445 1150}%
    \special{pa 1775 525}%
    \special{fp}%
    \special{pa 1445 1150}%
    \special{pa 525 1775}%
    \special{fp}%
    \special{pa 855 1150}%
    \special{pa 525 525}%
    \special{fp}%
    \special{pa 855 1150}%
    \special{pa 1775 1775}%
    \special{fp}%
    \special{pa 855 1150}%
    \special{pa 1445 1150}%
    \special{fp}%
    \special{pa 1150 265}%
    \special{pa 1150 118}%
    \special{fp}%
    \special{pa 2035 1150}%
    \special{pa 2182 1150}%
    \special{fp}%
    \special{pa 1150 2035}%
    \special{pa 1150 2182}%
    \special{fp}%
    \special{pa 265 1150}%
    \special{pa 118 1150}%
    \special{fp}%
    \graphtemp=.5ex\advance\graphtemp by 0.754in
    \rlap{\kern 0.312in\lower\graphtemp\hbox to 0pt{\hss $a$\hss}}%
    \graphtemp=.5ex\advance\graphtemp by 0.312in
    \rlap{\kern 0.754in\lower\graphtemp\hbox to 0pt{\hss $a$\hss}}%
    \graphtemp=.5ex\advance\graphtemp by 0.312in
    \rlap{\kern 1.546in\lower\graphtemp\hbox to 0pt{\hss $a$\hss}}%
    \graphtemp=.5ex\advance\graphtemp by 0.754in
    \rlap{\kern 1.988in\lower\graphtemp\hbox to 0pt{\hss $a$\hss}}%
    \graphtemp=.5ex\advance\graphtemp by 1.546in
    \rlap{\kern 1.988in\lower\graphtemp\hbox to 0pt{\hss $a$\hss}}%
    \graphtemp=.5ex\advance\graphtemp by 1.988in
    \rlap{\kern 1.546in\lower\graphtemp\hbox to 0pt{\hss $a$\hss}}%
    \graphtemp=.5ex\advance\graphtemp by 1.988in
    \rlap{\kern 0.754in\lower\graphtemp\hbox to 0pt{\hss $a$\hss}}%
    \graphtemp=.5ex\advance\graphtemp by 1.546in
    \rlap{\kern 0.312in\lower\graphtemp\hbox to 0pt{\hss $a$\hss}}%
    \graphtemp=.5ex\advance\graphtemp by 1.150in
    \rlap{\kern 0.000in\lower\graphtemp\hbox to 0pt{\hss $c'$\hss}}%
    \graphtemp=.5ex\advance\graphtemp by 0.754in
    \rlap{\kern 0.773in\lower\graphtemp\hbox to 0pt{\hss $b$\hss}}%
    \graphtemp=.5ex\advance\graphtemp by 0.000in
    \rlap{\kern 1.150in\lower\graphtemp\hbox to 0pt{\hss $c$\hss}}%
    \graphtemp=.5ex\advance\graphtemp by 0.754in
    \rlap{\kern 1.527in\lower\graphtemp\hbox to 0pt{\hss $b$\hss}}%
    \graphtemp=.5ex\advance\graphtemp by 1.150in
    \rlap{\kern 2.300in\lower\graphtemp\hbox to 0pt{\hss $c'$\hss}}%
    \graphtemp=.5ex\advance\graphtemp by 1.702in
    \rlap{\kern 1.462in\lower\graphtemp\hbox to 0pt{\hss $b$\hss}}%
    \graphtemp=.5ex\advance\graphtemp by 2.300in
    \rlap{\kern 1.150in\lower\graphtemp\hbox to 0pt{\hss $c$\hss}}%
    \graphtemp=.5ex\advance\graphtemp by 1.702in
    \rlap{\kern 0.838in\lower\graphtemp\hbox to 0pt{\hss $b$\hss}}%
    \graphtemp=.5ex\advance\graphtemp by 1.032in
    \rlap{\kern 1.150in\lower\graphtemp\hbox to 0pt{\hss $b$\hss}}%
    \graphtemp=.5ex\advance\graphtemp by 1.150in
    \rlap{\kern 0.855in\lower\graphtemp\hbox to 0pt{\hss $\bullet$\hss}}%
    \graphtemp=.5ex\advance\graphtemp by 1.150in
    \rlap{\kern 1.445in\lower\graphtemp\hbox to 0pt{\hss $\bullet$\hss}}%
    \hbox{\vrule depth2.300in width0pt height 0pt}%
    \kern 2.300in
  }%
}%
}
\caption{
The Petersen graph.\label{Pete}
}
\end{figure}

In this paper, we use the term \textit{double-star} for a $6$-vertex tree
having two adjacent vertices of degree $3$.  This is slightly non-standard; the
term is often used for any tree having two non-leaves.   Each edge of a
double-star incident to a leaf is a \textit{leaf edge}; the other edge is the
\textit{center edge}.  Given a subgraph $F$ whose components are pairwise
disjoint double-stars in a graph $G$, we say that $F$ {\em consists of
independent double-stars} if the set of leaves of components of $F$ is an
independent set of vertices in $G$.

In Figure~\ref{Pete}, we draw the Petersen graph on the torus with one
crossing.  The edge $c$ wraps top-to-bottom, and $c'$ wraps right-to-left.
The graph decomposes into a matching ($c$ and $c'$), an $8$-cycle ($a$), and a
double-star ($b$).  In view of Theorem~\ref{mainthm}, this decomposition proves
that the Petersen graph is $(7,2)$-edge-choosable.

\begin{definition}
A \emph{MED decomposition} of a $3$-regular graph is a decomposition of it
into subgraphs $G_1,G_2,G_3$, where $G_1$ is a Matching, the components of
$G_2$ are Even cycles, and $G_3$ consists of independent Double-stars.  Given a
MED decomposition, let $H$ be the graph formed from $G_3$ by deleting the
center edge from each component of $G_3$.
\end{definition}

In a MED decomposition, no vertex lies in both $G_1$ and $G_3$, since the third
edge at that vertex could not belong to any $G_i$.  Also, every vertex not in
$G_2$ belongs to the center edge of a component of $G_3$.  Finally, each
component of $H$ is a $3$-vertex path. 

To simplify the proof of Theorem~\ref{mainthm}, we prove two lemmas.  The first
establishes the sufficiency of an intermediate coloring by showing how to
complete the job of choosing a $2$-set $L$-coloring.  A {\it partial $2$-set
$L$-coloring} is a selection of at most $2$ colors from the list at each edge
so that the sets chosen at incident edges are disjoint.

\begin{lemma}\label{MED}
Let $G$ be a graph with maximum degree 3 that has a MED decomposition.
Given a $7$-uniform list assignment $L$ on $E(G)$, let $\phi$ be a partial
$2$-set $L$-coloring that chooses two colors for each edge in $G_1$ and one
color for each edge in $H$.  Each edge of $G_2$ on which $\phi$ does not save a
color is {\em needy}.  Each edge of $H$ incident to a needy edge is a
{\em sponsor}.  If the needy edges form a matching, the sponsors form a
matching, and each center edge incident to two sponsors has at least four
remaining colors, then $\phi$ extends to a $2$-set $L$-coloring.
\end{lemma}

\begin{proof}
We will choose one additional color for each edge of $H$ and two colors for
each center edge in $G_3$.  We must choose these colors to save one color on
each needy edge.  If this is done, then Corollary~\ref{even cycles} allows us
to choose colors for the edges of $G_2$ to complete the extension of $\phi$ to
a $2$-set $L$-coloring of $E(G)$.  For each edge $e\in E(G)$, let $L'(e)$
denote what remains of $L(e)$ when the colors chosen by $\phi$ on edges
incident to $e$ are deleted; this assignment $L'$ remains fixed throughout
the proof.

Let $F$ be a double-star in $G_3$ with center edge $e$, leaf edges $f_1$ and
$f_2$ at one end of $e$, and leaf edges $f_3$ and $f_4$ at the other end.  
We have $|L'(e)|\geq 3$.  Since $G_1$ and $G_3$ share no vertices, we have
$|L'(f_i)|\geq 5$ for all $i$.  Since the sponsors form a matching, there is at
most one sponsor in $\{f_1,f_2\}$ and at most one in $\{f_3,f_4\}$.

Suppose first that one leaf edge of $F$ is a sponsor.  We may assume by
symmetry that it is $f_1$; let $g$ be its incident needy edge.  Since $G_3$
consists of independent double-stars, the other endpoint of $g$ is
incident to an edge of $G_1$.  Thus, we have already chosen three colors for
edges incident to $g$.  By discarding additional colors if necessary, we may
assume that $|L'(g)|=4$.  Since $|L'(f_1)|\geq 5$, we save a color on $g$
by choosing a color in $L'(f_1)-L'(g)$ for $f_1$.  We choose
any two remaining colors for $e$ (at least two remain) and any remaining color
for each of $f_2$, $f_3$, and $f_4$.  Since $|L'(f_i)|\ge 5$, such choices are
available.

Suppose instead that we need to save a color both on $g$ incident to $f_1$
and on $g'$ incident to $f_3$.  As above, we have $|L'(f_1)|\geq 5$ and may
assume that $|L'(g)|=4$, and we can choose the second color for $f_1$ to save a
color on $g$.  Similarly, we can choose the second color for $f_3$ to save on
$g'$.  By hypothesis, $|L'(e)|\geq 4$ (before choosing for $f_1$ and $f_3$);
hence $e$ still has at least two available colors, and we choose two of them.
Finally, since $|L'(f_i)|\ge 4$, we can choose a second available color for
each of $f_2$ and $f_4$.

Now it remains only to choose two colors for each edge of $G_2$.  Since
completing the choices on $G_3$ produced a partial $2$-set $L$-coloring that
saves a color on each edge of $G_2$, these edges retain lists of size at least
4, and Corollary~\ref{even cycles} completes the coloring.
\end{proof}

The next lemma describes a step that may be repeated often in obtaining a
partial coloring to use as input to the procedure in Lemma~\ref{MED}.
The proof is similar to that of Lemma~\ref{key lemma}.

\begin{lemma}\label{cycle}
Given a $3$-regular graph $G$ with a MED decomposition, let $C$ be a cycle in
$G_2$ with edges $\VEC a1k$ in order.  Let $b_i$ be the other edge incident to
$a_i$ and $a_{i+1}$ (viewing indices modulo $k$).  Let $L$ be a list assignment
such that for each $i$, it holds that $|L(b_i)|\ge|L(a_i)|$ or $b_i\in E(H)$.
If $L(b_1)\not\esub L(a_1)$, or if $b_1\in E(H)$ and $|L(b_1)|<|L(a_1)|$, then
a color can be chosen for each edge in $\VEC b1k$ so that for each $i$, either
(1) a color is saved on $a_i$ or (2) $b_i\in E(H)$ and either
$|L(b_i)|<|L(a_i)|$ or some earlier $b_j$ is incident to $b_i$ in $H$.
\end{lemma}

\begin{proof}
If $L(b_1)\not\esub L(a_1)$, then choose $c\in L(b_1)-L(a_1)$ for $b_1$;
this saves a color on $a_1$.  Otherwise, $b_1\in H$ and $|L(b_1)|<|L(a_1)|$,
and we choose any $c\in L(b_1)$ for $b_1$ without saving a color on $a_1$.
In either case, remove $c$ from $L(a_2)$; if $c\notin L(a_2)$, then remove any
color from $L(a_2)$.  Also remove $c$ from the lists on edges incident to the
other end of $b_1$; these edges may be on $C$, on another even cycle in $G_2$,
or consist of another edge of $H$ and a central edge in a component of $G_3$.

Let $L'(e)$ denote the current remaining list at $e$.  In contrast to the proof
of Lemma~\ref{MED}, here we update $L'$ as we go along.  We proceed inductively
for increasing $i$.  After choosing a color for $b_{i-1}$, we have
$b_i\in E(H)$ or $|L'(b_i)|>|L'(a_i)|$, since we deleted a color from the list
at $a_i$ when choosing a color for $b_{i-1}$.

In the latter case, we choose a color for $b_i$ from $L'(b_i)-L'(a_i)$ and save
on $a_i$.  This includes the possibility that $b_i$ was visited earlier as
$b_j$, because in that case we reduced the available lists at both $b_j$ and
$a_i$ by choosing a color for $b_j$.  It also includes the possibility that
$b_i\in E(H)$ and $|L(b_i)|=|L(a_i)|$, which will be helpful later.

The remaining case is $b_i\in E(H)$ and $|L(b_i)|<|L(a_i)|$; the color deleted
from the list at $a_i$ does not guarantee $L'(b_i)-L'(a_i)\ne\emptyset$.
In this case, choosing a color from $L'(b_i)$ for $b_i$ may not save on $a_i$.

After choosing a color for $b_k$ by these rules, the proof is complete.
\end{proof}

In the case of Lemma~\ref{cycle} where choosing a color on $b_i$ does not save
a color on $a_i$, the edge $a_i$ will be designated needy, and the edge $b_i$
will be its sponsor.  In applying the lemma, we will sometimes exercise more
care in choosing the color for $b_i$ in that case.

We will use Lemma~\ref{cycle} to reach a state where Lemma~\ref{MED} applies.
We would like to save a color on each edge of $G_2$ while
choosing a partial $2$-set $L$-coloring with two colors chosen on each of the
other edges, as in Lemma~\ref{3ec}.  However, the double-stars outside these
even cycles cause complications and force us to leave needy edges and use
Lemma~\ref{MED}.

\begin{theorem}
\label{mainthm}
Every 3-regular graph $G$ having a MED decomposition is $(7,2)$-edge-choosable.
\end{theorem}
\begin{proof}
Before applying Lemma~\ref{MED}, we must save one color each on all edges of
$G_2$ except a matching, while confining the sponsors also to a matching and
choosing two colors for each edge of $G_1$ and one color for each edge of $H$. 
We do this in each component of $G_1\cup G_2\cup H$ independently; the central
edges of $G_3$ tie the graph together in the step performed later by
Lemma~\ref{MED}.  Let $J$ be a component of $G_1\cup G_2\cup H$.

\smallskip
{\bf Case 1.}
{\it $J$ does not contain two incident edges, one in $G_2$ and one in
$G_1\cup H$, that have distinct lists.}  In this case, for each component $R$
of $J\cap (G_1\cup G_2)$, the lists are identical on all edges of $R$ and all
edges of $H$ incident to it.  (There may be more than one such component $R$;
they may be joined by components of $H$.)  Choose any two colors $c$ and $c'$
from this list, and choose $\{c,c'\}$ for each edge of $R\cap G_1$.  For each
edge of $H$ incident to $V(R)$, choose $c$ or $c'$, avoiding the color on the
neighboring edge of $H$ if it has already been chosen by another such component.

Since $G_1$ and $H$ share no vertices, these choices are available.  Each edge
$e$ of $G_2$ is incident to two edges of $G_1\cup H$, at least one in $G_1$.
We have made at least three color assignments on edges incident to $e$ using
only two colors, thus saving at least one color on $e$.

\smallskip
{\bf Case 2.}
{\it $J$ has incident edges $e_1$ in $G_2$ and $e_2$ in $G_1\cup H$ with
distinct lists.} Let $C$ denote the even cycle in $J\cap G_2$ that contains
$e_1$.  Let $\VEC a1k$ be the edges of $C$ in order, with $a_1=e_1$ and $a_2$
incident to $e_1$.  Let $b_i$ be the edge incident to $a_i$ and $a_{i+1}$, so
$b_1=e_2$ and  $L(b_1)\neq L(a_1)$.  Note that $b_i\in E(G_1\cup H)$.  An edge
may be labeled as both $b_i$ and $b_j$ with $i\neq j$; this holds for each edge
of $G_1$ that is a chord of $C$.

We have chosen no colors on edges incident to $J$, so all lists have size $7$.  
Since $L(b_1)\ne L(a_1)$, Lemma~\ref{cycle} applies. 

{\it Step 1.  Choosing colors for $\VEC b1k$.}
Using Lemma~\ref{cycle}, we choose a color for each $b_i$ that saves on $a_i$
unless $b_i$ is incident to an edge of $H$ for which we have already chosen
a color.  That edge may for example be $b_j$ with $j<i$.  When we reach $b_i$,
the remaining list at $b_i$ may be smaller than that at $a_i$.

In this case, $a_i$ remains needy and $b_i$ is its sponsor.  This only happens
when $b_i$ is the second edge in its component of $H$ for which we are choosing
a color.  Let $e$ be the center edge incident to $b_i$.  If $e$ has six
remaining available colors when we choose for $b_i$, meaning that we have
not chosen a color for any edge incident to $e$ other than the other edge of
$H$ incident to $b_i$, then we choose any color for $b_i$ among its six
available colors.

If at most five colors remain available at $e$, then we choose a color for
$b_i$ to save a color on $e$.  This ensures that each center edge incident to
two sponsors will have at least four remaining colors, as needed for
Lemma~\ref{MED}.  With this refinement, we choose colors as in
Lemma~\ref{cycle}.

Let $L'$ denote the assignment of remaining available colors after this step,
and let $B=\{\VEC b1k\}$.  For any edge $e$ incident to $B$, each color $c$
chosen on an incident edge $e'$ of $B$ eliminates one element from $L'(e)$
(which can be any element if $c$ was not in the list); it also eliminates $c$
from $L'(e')$.  If $e\in E(G_2)$, then $e$ may be incident to another edge of
$B$.  Thus the list for $e$ decreases by at least as much as the list for $e'$.

{\it Step 2.  Handling additional components of $G_2$ in $J$.}
If we have not yet chosen colors on all of $J$, then $J\cap G_2$ contains
another unprocessed cycle $C'$ with length $k'$, reachable from an already
processed cycle $C$ via an edge of $G_1$ or via a component (two edges) of $H$
(or possibly reachable in both ways).  Suppose that $\VEC b1k$ were the
edges incident to $C$; we chose a color for each of these while processing $C$.

If some $b_i$ is incident to $C'$ such that $L'(b_i)$ has an element outside
the list of some edge in $C'$ incident to $b_i$, then let $b'_1=b_i$; in this
case $b'_1\in E(G_1)$.  Alternatively, if $f_1$ and $f_2$ are incident
edges in $H$, with $f_1=b_i$ and $f_2$ incident to $C'$, then let $b'_1=f_2$.
In either case, index the edges of $C'$ as $\VEC {a'}1{k'}$ in order, starting
with edges incident to $b'_1$ as $a'_1$ and $a'_2$, choosing $a'_1$ so that
$L'(b'_1)\not\esub L'(a'_1)$ if possible.

We selected $b'_1\in E(G_1)$ only when $L'(b'_1)\not\esub L'(a'_1)$.
Also, previous applications of Lemma~\ref{cycle} have left
$|L'(b'_i)|\ge |L'(a'_i)|$ unless $b'_i\in E(H)$ and the edge incident to $b'_i$
in $H$ already had a color chosen while processing an earlier cycle.  This 
requires $b'_1\in E(H)$ and $|L'(b'_1)|<|L'(a'_1)|$.  Hence in either case the
hypotheses for Lemma~\ref{cycle} hold.  We return to Step 1 to process $C'$,
treating $C',k',b'_i,a'_i$ as $C,k,b_i,a_i$.  

Continuing in this way, we process all components of $J\cup G_2$ unless we
reach a state where a remaining cycle $C$ in $J\cup G_2$ is incident only to
edges of $G_1$, and the remaining list on each such edge is the same as the
remaining list on the incident edges of $C$.  Since each vertex of $C$ has
an incident edge in $G_1$, all the lists on $C$ and its incident edges are the
same.  

That list contains six colors, since we have only chosen one color on each incident edge in $G_1$.  Since each edge of $C$ has six remaining colors in its list, and we have already chosen one color on each incident edge in $G_1$,
we have
already saved a color on each edge of $C$.  We choose an additional color from
this list as the second color for all the edges of $G_1$ incident to $C$.

Under the algorithm used in Lemma~\ref{cycle}, only the second edge of a
component of $H$ on which a color is chosen can become a sponsor.  Hence the
sponsors form a matching.  Furthermore, an edge $a_i$ of $G_2$ becomes needy
only when the edge incident to $a_i$ at its endpoint incident to the subsequent
edge ($a_{i+1}$ when processing its cycle in $G_2$) is a sponsor.  This
prevents two incident edges of $G_2$ from being left needy, since $G_3$
consists of independent double-stars.  We have therefore produced a partial
$2$-set $L$-coloring that satisfies the hypotheses of Lemma~\ref{MED}, and
Lemma~\ref{MED} completes the $2$-set $L$-coloring.
\end{proof}

We have used Theorem~\ref{mainthm} to check that all snarks with at most 24
vertices are $(7,2)$-edge-choosable, as are the double star snark, the Szekeres
snark, the Goldberg snark, the Watkins snarks of orders 42 and 50, all
cyclically 5-edge-connected snarks of order 26, and the infinite family of
flower snarks.  Cavicchioli et al.~\cite{cav1} presented drawings of the 67
pairwise non-isomorphic snarks with at most 24 vertices as well as all
cyclically 5-edge-connected snarks of order 26.  Their drawings illustrate that
the snarks have cycles omitting two adjacent vertices.  Since a snark has girth
at least five, the leaves of any single double-star form an independent set.
Thus every snark having a cycle of length $n-2$ has a MED decomposition using
this cycle.  All snarks listed above have such decompositions.

\begin{example}
\rm
Isaacs~\cite{isaacs} discovered the flower snarks.  For odd $k$ with $k\ge5$,
the flower snark $F_k$ consists of $4k$ vertices, labeled $w_i,x_i,y_i,z_i$ for
$1\le i\le k$.  The edge set consists of a cycle $[\VEC z1k]$ of length $k$, a
cycle $[\VEC x1k,\VEC y1k]$ of length $2k$, and the $k$ stars with center $w_i$
and leaves $x_i,y_i,z_i$.  To prove that $F_k$ is (7,2)-edge-choosable, we
provide a MED decomposition; as above, it suffices to find a cycle of length
$4k-2$ that omits adjacent vertices.  For example,
%
$[x_1,x_2,w_2,y_2,y_3,w_3,x_3,\ldots,x_{k-1},w_{k-1},y_{k-1},y_k,w_k,z_k,\ldots,z_1,w_1]$
omits only $x_k$ and $y_1$.
\end{example}

We do not know whether every snark has a cycle that omits just two adjacent vertices.  
We conjecture at least that all snarks have MED decompositions.  
In fact, snarks seem to be easy graphs in which to find MED decompositions, and we believe the following.

\begin{conjecture}
\label{mainconj}
Every 2-connected graph with maximum degree 3 has a MED decomposition (a
decomposition into a matching, a union of disjoint even cycles, and an
independent set of double-stars).
\end{conjecture}

If Conjecture~\ref{mainconj} is true, then by Theorem~\ref{mainthm}
all 2-connected graphs with maximum degree 3 are $(7,2)$-edge-choosable.
The argument does not extend to graphs with cut-vertices.



\begin{example}\label{gstar}
\rm
{\it There are $3$-regular graphs having no MED decomposition.}
A small example is the standard small $3$-regular graph having no perfect
matching.  Let $H$ be the graph obtained from $K_4$ by subdividing one edge.
Form $G^*$ from $3H$ by adding one vertex $v$ made adjacent to the $2$-valent
vertices of $3H$.  Let $x_1,x_2,x_3$ be the neighbors of $v$.

Suppose that $G^*$ has a MED decomposition.  Putting an $x_i$ on an even cycle
leaves a claw, which cannot be decomposed.  Hence the decomposition must
have each $x_i$ as a central vertex in a double-star.  Since a triangle can't
be decomposed, the other central vertex must be $v$.  Now the three required 
double-stars share $v$, yielding a contradiction.
\end{example}

Nevertheless, $G^*$ is $(7,2)$-edge-choosable and hence is not a counterexample
to Mohar's conjecture.  The proof unfortunately requires some case analysis.

\begin{theorem}
The graph $G^*$ is $(7,2)$-edge-choosable.
\end{theorem}
\begin{proof}
Again we reduce the problem to $(4,2)$-edge-choosability of even cycles.

Let $H'$ be the graph obtained from $K_4+K_1$ by subdividing one edge and
adding a pendant edge joining the new vertex to the isolated vertex.  Let $a$
be the pendant edge, incident to edges $b_1$ and $b_2$.  Let $d$ be the edge not
incident to $b_1$ or $b_2$.  Let $c_1$, $c_2$, $c_3$, $c_4$ be the edges of
the remaining cycle $C$ in order, with $c_1$ and $c_2$ incident to $b_1$.

Given a $7$-uniform list assignment $L$ on the edges of $H'$, a partial $2$-set
$L$-coloring can be chosen with one color each on $b_1$ and $b_2$ and two
on $d$ so that a color is saved on each edge of $C$.  If all the lists except
that on $a$ are the same, then we pick colors for $b_1$ and $b_2$ and use the 
same two colors
on $d$.  Otherwise, since all the lists have size $7$, the argument of
Lemma~\ref{cycle} finds such a coloring.  After doing this, only one remaining
color at $b_2$ is forbidden (by $b_1$), so we have at least five choices for a
second color there. 

Now let $H'_1, H'_2, H'_3$ denote the three copies of $H'$ in $G^*$, and let
$a_i$ denote edge $a$ in $H'_i$.  First, choose colors on $H'_1$ as specified by
the claim above.  When choosing colors on $H'_2$ in this way, we pause before
making the second choice on $b_2$.  There are at least five choices for this
color, and at least four choices remain at $a_1$, so we can choose color
$\alpha$ for $a_1$ and $\beta$ for $b_2$ to save a color on $a_2$.

Let $L'(e)$ denote the list of available colors for $e$ at this time.  Since we
have chosen a color on $a_1$ and three colors incident to $a_2$ in $H'_2$, but
saved a color on $a_2$, we may assume that $|L'(a_2)|=4$.  Similarly, from the
colors chosen on $H'_1$ we may assume that $|L'(a_1)|=3$.

If we can choose a $2$-set coloring for $E(H_3)$ so that $a_3$ receives a color
not in $L'(a_2)$, then we can safely pick one more color on $a_1$ and two on
$a_2$.  We can then return to pick the last extra color needed on $H'_1$ and
$H'_2$ outside the copies of $C$, after which we can color each copy of $C$
since we save a color on each of its edges.

Hence it remains to choose the colors on $H_3$ so that $a_3$ has a color
outside a specified set of size $5$ (avoiding also the color on $a_1$).
There remain two desirable colors in $L(a_3)$; call them $c$ and $c'$.
We show that $H'_3$ can be colored using $c$ or $c'$ on $a_3$.

Let $L_0(e)=L(e)-\{c,c'\}$ for $e\in E(H_3)$.  If $L_0(b_1)$ has a color
$\beta$ outside $L(c_1)\cap L(c_2)$, then our original claim about $H'$ using
the argument of Lemma~\ref{cycle} allows us to pick $\beta$ for $b_1$, two
colors for $d$, and two colors for $b_2$ that are not both in $\{c,c'\}$ (while
saving on each edge of $C$).  Now we can use $c$ or $c'$ and one of three other
colors on $a_3$ and then return to complete the coloring on $H_3$.

The same outcome occurs if $L_0(d)$ has a color not in $L(c_1)\cap L(c_2)$.
Hence $L(c_1)\cap L(c_2)$ contains $L_0(b_1)\cup L_0(d)$.  The restricted sets
have size at least $5$, and $|L(c_1)\cap L(c_2)|\le 7$, so
$|L_0(b_1)\cap L_0(d)|\ge3$.  The symmetric argument switching the roles of
$b_1$ and $b_2$ shows that $|L_0(b_2)\cap L_0(d)|\ge3$.

Hence we may choose $\alpha\in L_0(b_1)\cap L_0(d)$ and
$\beta\in L_0(b_2)\cap L_0(d)-\{\alpha\}$, putting color $\alpha$ on $b_1$,
$\beta$ on $b_2$, and both on $d$.  Choose another color on $b_2$, then $c$ or
$c'$ on $a_3$, then a color that remains for $b_1$.  This leaves lists of size
$4$ on $C$, and we can complete the $2$-set $L$-coloring.
\end{proof}

Mohar's conjecture remains open, and Example~\ref{gstar} shows that something
stronger than Theorem~\ref{mainthm} will be needed to prove it.  There may be a
generalization of MED decomposition that suffices for
$(7,2)$-edge-choosability, permitting more complicated structures than
double-stars in addition to the matching and even cycles.
Whether this is true or not, it would still be interesting to show that 
every 2-connected graph with maximum degree 3 has a MED decomposition.


\end{document}